\documentclass[11pt]{article}
\usepackage{amsfonts}
\usepackage{latexsym}
\usepackage{amsmath}
\usepackage{amssymb}
\usepackage{amsthm}
\newtheorem{theorem}{Theorem}[section]

\newtheorem{proposition}[theorem]{Proposition}
\newtheorem{corollary}[theorem]{Corollary}

\theoremstyle{definition}
\newtheorem{definition}[theorem]{Definition}

\theoremstyle{remark}

\newtheorem*{note*}{Note}



\begin{document}

\title{\bf 
R\'enyi Divergence and $L_p$-affine surface area  for convex bodies
}

\author{ Elisabeth M. Werner 
\thanks{Partially supported by an NSF grant, a FRG-NSF grant and  a BSF grant}}

\date{}

\maketitle

\begin{abstract}
We show that the fundamental  objects of the $L_p$-Brunn-Minkowski theory, 
namely the $L_p$-affine surface areas for a convex body,   are 
closely related to information theory: 
they are exponentials of R\'enyi divergences  of the cone measures of  a convex body and its polar.
\par
We give geometric interpretations for  all R\'enyi divergences $D_\alpha$,  
not just for the previously treated special case of relative entropy which is the case $\alpha =1$.
Now,  no symmetry assumptions are needed and, if at all, only very weak regularity assumptions are required.
\par
Previously, the relative entropies appeared only after performing   second order expansions of certain expressions. Now   already  first order expansions makes them appear. Thus,  in the new approach we detect ``faster"  details  about the  boundary of a convex body.

\end{abstract}

      \vskip 2cm \noindent

\section{Introduction.}
There exists a fascinating connection between convex geometric analysis and information
theory. An example is the close parallel between geometric inequalities for convex
bodies and inequalities for probability densities. For instance, the Brunn-Minkowski
inequality and the entropy power inequality follow both in a very similar way from the
sharp Young inequality (see. e.g., \cite{DCT}).
\par
In several recent papers, Lutwak, Yang, and Zhang \cite{LYZ2000, LYZ2002/1, LYZ2004/1, LYZ2005} established further connections
between convexity and information theory. For example, they showed in \cite{LYZ2002/1}  that the
Cramer-Rao inequality corresponds to an inclusion of the Legendre ellipsoid and the
polar $L_2$-projection body. The latter is a basic notion from the  $L_p$-Brunn-Minkowski
theory. 
This $L_p$-Brunn-Minkowski theory has its origins in the 1960s when Firey introduced
his $L_p$-addition of convex bodies. It 
evolved rapidly over the last years and
due to a number of highly influential works (see, e.g., 
\cite{Ga3}, \cite{GaZ} - \cite{HLYZ}, 
\cite{Klain1}, \cite{Klain2},
\cite{Lud2} - \cite{Lu1988},  \cite{LYZ2002}, \cite{LYZ2004}, 
\cite{LZ} -
\cite{NPRZ}, \cite{RuZ}, 
\cite{Schu} - 
\cite{WY2008},
\cite{Z3}),  is now  a central part of modern convex geometry. In fact,   this theory
redirected much of the research about  convex bodies from  the Euclidean aspects 
 to the study of the affine geometry of these bodies,  and some questions that had been considered Euclidean in nature turned out to be affine problems. For example,  the famous Busemann-Petty Problem (finally laid to rest in
\cite{Ga1, GaKoSch,  Z1, Z2}), was shown to be an affine problem with the introduction of intersection bodies by Lutwak in \cite{Lu1988}.
\par
Two fundamental notions within the
$L_p$-Brunn-Minkowski theory are $L_p$-affine surface areas,   introduced by Lutwak   
in the ground
breaking paper \cite{Lu2} and $L_p$-centroid bodies introduced by Lutwak and Zhang in \cite{LZ}. See Section 3 for the definition of those quantities.
\par
Based on these quantities, 
Paouris and Werner  \cite{PW1} established yet another  relation between affine convex
geometry and information theory. They proved that the exponential of the
relative entropy of the cone measure of a symmetric convex body and its polar equals a limit of
normalized $L_p$-affine surface areas. 
Moreover, also in  \cite{PW1}, 
Paouris and Werner gave  geometric interpretations of the relative entropy of the cone measures of a {\em sufficiently smooth, symmetric } convex body and its polar.
\par
In this paper we show  that the very core of the $L_p$-Brunn-Minkowski theory, namely the $L_p$-affine surface areas itself,  are 
concepts of information theory: 
They are exponentials of R\'enyi divergences  of the cone measures of a convex body  and its polar.
This identification allows to translate known properties from one theory to the other. 
\par
Even more is  gained.  
Geometric interpretations for {\em all}  R\'enyi divergences $D_\alpha$ of  cone measures of a convex body  and its polar are given for {\em all}  $\alpha$, not just for the special case of relative entropy which corresponds to the case $\alpha =1$.
We refer to Sections 2 and 3 for the definition of $D_\alpha$. 
No symmetry assumptions on $K$ are needed. Nor do these new geometric interpretations  require the strong smoothness assumptions
of \cite{PW1}.
\par
In the context of the  $L_p$-centroid bodies, the relative entropies appeared only after performing   {\em second order expansions} of certain expressions. The remarkable fact now  is that in our approach here, already  {\em first order expansions} makes them appear. Thus,  these bodies detect ``faster"  details  of the  boundary of a convex body than the  $L_p$-centroid bodies.

\vskip 4mm
The paper is organized as follows. 
In Section 2 we introduce R\'enyi divergences for convex bodies and describe some of their properties.
We also introduce $L_p$-affine surface areas and mixed  $p$-affine surface areas. 
\par
The main observations are Theorems \ref{p-aff=renyi}  and \ref{mixed-aff=renyi} which show
that $L_p$-affine surface areas and mixed  $p$-affine surface areas are exponentials of R\'enyi divergences.
These identifications allow to translate known properties from one theory to the other - this is done in the rest of Section 2 and in Section 3. 
Also, in Section 3, we give geometric interpretations for R\'enyi divergences $D_\alpha$ of cone measure of convex bodies  for all $\alpha$, including new ones for the 
relative entropy not requiring  the (previously necessary)  strong smoothness and symmetry assumptions on the body.

\vskip 4mm
{\bf Further Notation.}
\vskip 2mm
\noindent
Throughout the paper, we will  assume that  the centroid of a
convex body $K$ in $\mathbb R^n$ is at the origin. 
We work in ${\mathbb R}^n$, which is
equipped with a Euclidean structure $\langle\cdot ,\cdot\rangle $.
We denote by $\|\cdot \|_2$ the corresponding Euclidean norm.  $B^n_2(x, r)$ is the ball centered at $x$ with radius $r$. We
write $B_2^n=B^n_2(0,1)$ for the Euclidean unit ball centered at $0$ and $S^{n-1}$ for the
unit sphere. Volume is denoted by $|\cdot |$ or, if we want to emphasize the dimension,  by $\text{vol}_d(A)$ for a $d$-dimensional set $A$. 
$K^\circ=\{ y\in \mathbb{R}^n: \langle x, y \rangle \leq 1 \  \text{for all }  \  x \in K\}$ is the polar body of $K$.
\par
For a point $x \in \partial K$, the boundary of $K$,  $N_K(x)$ is the outer unit normal 
in $x$ to $K$ and $\kappa_K(x)$  is the (generalized) Gauss curvature in $x$.
We  write $K\in
C^2_+$,  if $K $ has $C^2$ boundary $\partial K$  with everywhere strictly
positive Gaussian curvature $\kappa_K$.  $\mu_K$ is the usual surface area measure on $\partial K$.
 $\sigma $ is the usual surface area measure on $S^{n-1}$.
\par
Let  $K$ be a convex body in $\mathbb{R}^n$ and let  $u \in S^{n-1}$. Then
$h_{K}(u)$ is the support function of direction $u\in S^{n-1}$,
and $f_{K}(u)$ is the curvature function, i.e. the reciprocal of
the Gaussian curvature $\kappa _K(x)$ at this point $x \in
\partial K$ that has $u$ as outer normal.

\vskip 5mm
\section{R\'enyi divergences for convex bodies.}
Let $(X, \mu)$ be a measure space  and let  $dP=pd\mu$ and  $dQ=qd\mu$ be probability measures on $X$ that are  absolutely continuous with respect to the measure $\mu$. 
Then the R\'enyi divergence of order $\alpha$,  introduced by  R\'enyi \cite{Ren}  for $\alpha >0$,  is defined as 
\begin{equation}\label{renyi}
D_\alpha(P\|Q)=
\frac{1}{\alpha -1} \log \int_X p^\alpha q^{1-\alpha} d\mu,
\end{equation}
It is  the convention to put  $p^\alpha q^{1-\alpha}=0$, if $p=q=0$, even if $\alpha <0$ and $\alpha >1$. 
The integrals 
\begin{equation}\label{HellInt}
\int_X p^\alpha q^{1-\alpha} d\mu
\end{equation}
are  also  called {\em Hellinger integrals}. See e.g.  \cite {Liese/Vajda2006} for those integrals  and  additional information.
\par 
Usually, in the literature, $\alpha \geq 0$. However, we will  also consider $\alpha <0$, provided the expressions exist.
We normalize the measures  as, again  usually in the literature, the measures  are probability measures. 
\vskip 4mm
{\bf Special cases.}
\vskip 2mm
\noindent
(i) The case $\alpha=1$ is also called the 
{\em Kullback-Leibler divergence} or {\em relative entropy} from $P$ to $Q$ (see \cite{CT2006}). It is  obtained as the limit as $\alpha \uparrow 1$ in  (\ref{renyi}) and one gets
\begin{equation}\label{relent}
 D_{KL}(P\|Q)= D_{1}(P\|Q)=  \lim _{\alpha \uparrow 1} D_\alpha(P\|Q)= \int_{X} p \log \frac{p}{q} d\mu.
\end{equation}
(The limit  $\alpha \rightarrow 1$ may not exist but  limit $\alpha \uparrow 1$ exists \cite{Liese/Vajda1987}).
\vskip 2mm
\noindent
(ii) The case $\alpha=0$ gives   for $q \neq 0$ (with the convention that $0^0=1$ ) that

\begin{equation}\label{0}
  D_{0}(P\|Q)=  0,
\end{equation}
as $dQ=q d\mu$ is a probability measure on $X$. If $q=0$, then $D_{0}(P\|Q)= - \infty$.
\vskip 2mm
\noindent
(iii) The case $\alpha=\frac{1}{2}$ gives
\begin{equation}\label{halb}
  D_{\frac{1}{2}}(P\|Q)=  D_{\frac{1}{2}}(Q\|P) = -2 \log \int_{X} p ^\frac{1}{2} q ^\frac{1}{2}   d\mu.
\end{equation}
The expression $\int_{X} p ^\frac{1}{2} q ^\frac{1}{2}   d\mu$ is also called the {\em Bhattcharyya  coefficient} or 
{\em Bhattcharyya distance} of $p$ and $q$.
\vskip 2mm
\noindent
(iii) The cases $\alpha=\infty$ and $\alpha= - \infty$. 
\begin{equation}\label{infty}
  D_{\infty}(P\|Q)=    \log  \  \left(\sup_x \text{ess} \  \frac{p(x)}{q(x)}\right),      
\end{equation}
and
\begin{equation}\label{minusinfty}
D_{-\infty}(P\|Q)=   -  \  \left(\sup_x \text{ess} \   \frac{q(x)}{p(x)}\right) = - D_{\infty}(Q\|P).
\end{equation}

\vskip 5mm
Note that  for all $- \infty \leq \alpha \leq \infty$, $\alpha \neq 1$, 
\begin{equation}\label{PgleichQ}
D_\alpha(Q\|P)= \frac{\alpha}{1-\alpha } D_{1-\alpha}(P\|Q).
\end{equation}
As $\alpha \uparrow 1$, the limit on the  left and the limit on the right of (\ref{PgleichQ}) exist  and are  equal and equal to  $ D_{1}(Q\|P)= \int_{X} q \log \frac{q}{p} d\mu$. 
Thus (\ref{PgleichQ}) holds for all $- \infty \leq \alpha \leq \infty$.

\vskip 5mm
\noindent
We will now  consider  R\'enyi divergence for  convex bodies $K$
 in $\mathbb R^n$. 
Let 
\begin{equation}\label{densities}
p_K(x)= \frac{ \kappa_{K}(x)}{\langle x, N_{K}(x) \rangle^{n} \  n|K^{\circ}|} \, , \   \ q_K(x)= \frac{\langle x, N_{K}(x) \rangle }{n\  |K|}.
\end{equation}
\noindent 
Then 
\begin{equation}\label{PQ}
P_K=p_K\  \mu_K \ \ \ \text{and}   \ \ \   Q_K=q_K \ \mu_K
\end{equation}
 are probability measures on $\partial K$ that are absolutely continuous with respect  to $\mu_{K}$.
\par
Recall  that the normalized cone measure $cm_K$
on $\partial K$ is defined as follows:
For every measurable set $A \subseteq \partial K$
\begin{equation}\label{def:conemeas} 
cm_{K}(A)  = \frac{1}{|K|} \bigg| \big\{ta : \ a \in A, t\in [0,1] \big\} \bigg|.
\end{equation}
\vskip 3mm
\noindent
The next proposition is well known. See e.g. \cite {PW1} for a proof. It shows that the measures $P_K$ and $Q_K$ defined in (\ref{PQ})
 are the cone measures  of $K$  and $K^\circ$.  $N_K:\partial K \rightarrow S^{n-1}$, $x \rightarrow N_K(x)$  is the Gauss map.
\vskip 3mm
\noindent 
\begin{proposition} \label{prop:conemeas}
\noindent 
Let $K$ a  convex body in $\mathbb R^n$. Let $P_K$ and $Q_K$ be the probability measures on $\partial K$  defined by (\ref{PQ}).  Then
$$
Q_K= cm_{K},
$$
\noindent 
or, equivalently, for every measurable subset $A$ in $ \partial K$
$Q_K(A)= cm_{ K}(A)$.
\newline
If $K$ is in addition in $C^2_+$, then 
$$
P_K= N_{K}^{-1}N_{K^{\circ}}cm_{K^{\circ}}
$$
\noindent 
or, equivalently, for every measurable subset $A$ in $ \partial K$
$$ 
P_K(A)= cm_{ K^{\circ}} \bigg(N_{{K^{\circ}}}^{-1} \big(N_{K} (A)\big)\bigg).
$$
\end{proposition}

\vskip 5mm 
For $\alpha =1$, the relative entropy  of a convex body $K$ in  $\mathbb{R}^n$ was considered in \cite{PW1}, namly
\begin{eqnarray*}\label{KL1}
D_{1}(P_K\|Q_K) &= & D_{KL}(P_K\|Q_K) \\ &=& \int_{\partial K}
\frac{\kappa_K(x)} {n |K^\circ| \langle x, N_{K}(x) \rangle ^n} \log \left(\frac{|K|\kappa_K(x)}{|K^\circ|\langle x, N_{K}(x) \rangle ^{n+1}}\right) d \mu_K(x)
\end{eqnarray*}
\begin{eqnarray*}\label{KL2}
D_{1}(Q_K\|P_K) &= & D_{KL}(Q_K\|P_K) \\ &= &\int_{\partial K} \frac{\langle x, N_{K}(x) \rangle}{n|K|} \log\left(\frac{|K^\circ|\langle x, N_{K}(x) \rangle ^{n+1}}{|K| \kappa_K(x)} \right)d \mu_K(x),
\end{eqnarray*}
provided the expressions exist. 
\vskip 5mm
 We now define 
 the R\'enyi divergence of $K$ of order $\alpha$ for all other $\alpha$,  $- \infty \leq  \alpha \leq  \infty$,  $\alpha \neq 1$.
 \begin{definition}\label{Renyidiv} 
{\em  Let $K$ be a convex body in $\mathbb{R}^n$ and let $- \infty < \alpha < \infty$,  $\alpha \neq 1$. Then the R\'enyi divergences of order $\alpha$ of $K$ are
\begin{equation}\label{K-renyi1}
D_\alpha(Q_K\|P_K)= \frac{1}{\alpha -1} \log\left( \frac{
\int _{\partial K} \frac{\kappa_K^{1-\alpha} d\mu_K}{\langle x, N_K(x) \rangle ^{n-\alpha(n+1)}}}{n |K|^\alpha |K^\circ|^{1-\alpha}} \right)
\end{equation}
\vskip 2mm
\noindent
\begin{equation}\label{K-renyi2}
D_\alpha(P_K\|Q_K)= \frac{1}{\alpha -1} \log\left( \frac{\int _{\partial K}
\frac{\kappa_K^{\alpha} d\mu_K}{\langle x, N_K(x) \rangle ^{\alpha(n+1) -1}}}{n |K|^{1-\alpha} |K^\circ|^\alpha} 
\right)
\end{equation}
\vskip 2mm
\noindent
\begin{equation}\label{K-renyi3}
D_\infty(Q_K\|P_K)= \log\left( \sup_{x \in \partial K} \text{ess}\    \frac{ |K^\circ| \  \langle x, N_K(x) \rangle ^{n+1}}{  |K|  \  \kappa_K(x)}\right)
\end{equation}
\vskip 2mm
\noindent
\begin{equation}\label{K-renyi4}
D_\infty(P_K\|Q_K)= \log\left( \sup_{x \in \partial K} ess\   \frac{  |K|  \  \kappa_K(x)} { |K^\circ| \  \langle x, N_K(x) \rangle ^{n+1}}\right)
\end{equation}
and 
\vskip 2mm
\noindent
\begin{equation}\label{K-renyi5}
D_{-\infty}(Q_K\|P_K)= - D_\infty(P_K\|Q_K), \hskip 1mm D_{-\infty}(P_K\|Q_K)= - D_\infty(Q_K\|P_K),
\end{equation}
provided the expressions exist.}
\end{definition}
\vskip 5mm
{\bf Remarks.} 
\vskip 2mm
\noindent
(i) By (\ref{PgleichQ}) for all $- \infty \leq \alpha \leq \infty$, $\alpha \neq 1$,
\begin{equation*}\label{gleich}
D_\alpha(Q_K\|P_K)= \frac{\alpha}{\alpha -1} D_{1-\alpha}(P_K\|Q_K).
\end{equation*}
This identity also holds for $\alpha \uparrow 1$.
Therefore, it is enough to consider only one of the two, $D_\alpha(Q_K\|P_K)$ or $D_\alpha(P_K\|Q_K)$.
\vskip 3mm
\noindent
(ii) If we  put  $N_K(x)=u \in S^{n-1}$,  then  $\langle x, N_{K}(x)\rangle = h_K(u)$. If $K$ is in $C^2_+$,   then  $d\mu_K= f_K d \sigma$. Hence, in that case, we can express the R\'enyi divergences also as
\begin{equation}
D_\alpha(Q_K\|P_K)= \frac{1}{\alpha -1} \log\left(\frac{\int _{ S^{n-1}} \frac{f_K(u)^\alpha d\sigma(u)}{h_K(u)  ^{n-\alpha(n+1)}}}{n |K|^\alpha |K^\circ|^{1-\alpha}}  \right)
\end{equation}
\vskip 2mm
\noindent
\begin{equation}
D_\alpha(P_K\|Q_K)= \frac{1}{\alpha -1} \log\left( \frac{\int _{ S^{n-1}}
\frac{f_K(u)^{1-\alpha} d\sigma(u)}{h_K(u)  ^{\alpha(n+1) -1}}
}{n |K|^{1-\alpha} |K^\circ|^\alpha} 
\right)
\end{equation}
\vskip 2mm
\noindent
Accordingly for $D_{KL}(Q_K\|P_K)$ and $D_{KL}(P_K\|Q_K)$.
\vskip 5mm
Let $K_1, \dots K_n$ be convex bodies in $\mathbb{R}^n$. Let $u \in S^{n-1}$. For $1 \leq i \leq n$, define
\begin{equation}\label{i-densities}
p_{K_i}(u)= \frac{ 1}{ n^\frac{1}{n} |K_i^{\circ}|^\frac{1}{n} h_{K_i}(u)} \, , \   \ q_{K_i}(u)= \frac{f_{K_i}(u)^\frac{1}{n} h_{K_i}(u)^\frac{1}{n} }{n^\frac{1}{n}\  |K_i|^\frac{1}{n}}.
\end{equation}
and measures on $S^{n-1}$ by
\begin{equation}\label{iPQ}
P_{K_i}=p_{K_i}\  \sigma \ \ \ \text{and}   \ \ \   Q_{K_i}=q_{K_i} \ \sigma.
\end{equation}
\vskip 4mm
\noindent
Then we  define the R\'enyi divergences of order $\alpha$ for convex bodies $K_1, \dots K_n$ by
\vskip 4mm
\noindent
\begin{definition}\label{mixed}
{\em \  
Let $K_1, \dots K_n$ be convex bodies in $\mathbb{R}^n$. Then for $-\infty < \alpha < \infty$, $\alpha \neq 1$
\begin{equation*}
D_\alpha(Q_{K_1} \times \dots \times Q_{K_n} \|P_{K_1} \times \dots \times P_{K_n})=  \frac{\log\left(  \int_{S^{n-1}} \prod_{i=1}^n \frac{f_{K_i}^{\frac{\alpha}{n} }h_{K_i}^{\frac{\alpha}{n}-(1-\alpha)}} { n^\frac{1}{n} \ |K_i|^\frac{\alpha}{n} |K_i^\circ|^\frac{1-\alpha}{n}} d\sigma\right)}{\alpha -1}
\end{equation*}
\vskip 2mm
\noindent
\begin{equation*}
D_\alpha(P_{K_1} \times \dots \times P_{K_n} \|Q_{K_1} \times \dots \times Q_{K_n})=  \frac{\log\left(  \int_{S^{n-1}} \prod_{i=1}^n \frac{f_{K_i}^{\frac{1-\alpha}{n} }h_{K_i}^{\frac{1-\alpha}{n}-\alpha}} { n^\frac{1}{n} \ |K_i|^\frac{1-\alpha}{n} |K_i^\circ|^\frac{\alpha}{n}} d\sigma\right)}{\alpha -1} 
\end{equation*}
provided the expressions exist. 
\newline
For $\alpha =1$ the definitions were given in \cite{PW1}:
\vskip 2mm
\noindent
\begin{eqnarray*}
&&\hskip -10mm D_1(Q_{K_1} \times \dots \times Q_{K_n} \|P_{K_1} \times \dots \times P_{K_n})=  \\ 
&&\hskip 19mm \int_{S^{n-1}} \prod_{i=1}^n \frac{f_{K_i}^{\frac{1}{n} }h_{K_i}^{\frac{1}{n}}} { n^\frac{1}{n} \ |K_i|^\frac{1}{n} }  \log\left(\prod_{i=1}^n \frac{ |K_i^\circ |^\frac{1}{n}\  f_{K_i}^{\frac{1}{n} }h_{K_i}^{1+\frac{1}{n}}} {  |K_i|^\frac{1}{n} }  \right) d\sigma
\end{eqnarray*}
\begin{eqnarray*}
&&\hskip -10mm D_1(P_{K_1} \times \dots \times P_{K_n} \|Q_{K_1} \times \dots \times Q_{K_n}) =  \\ 
&& \hskip 19mm \int_{S^{n-1}} \prod_{i=1}^n \frac{h_{K_i}^{-1} }{ n^\frac{1}{n} \ |K_i^\circ|^\frac{1}{n} }  \log\left(\prod_{i=1}^n \frac{ |K_i |^\frac{1}{n} } {|K_i^\circ |^\frac{1}{n}f_{K_i}^{\frac{1}{n}} h_{K_i}^{1+\frac{1}{n}} }  \right) d\sigma, 
\end{eqnarray*}
provided the expressions exist.
}
\end{definition}
\vskip 4mm
{\bf Remark.} 
For $- \infty < \alpha < \infty$,  $\alpha \neq 1$, 
\begin{eqnarray}\label{gleich2}
&& \hskip -10mm D_\alpha(P_{K_1} \times \dots \times P_{K_n} \|Q_{K_1} \times \dots \times Q_{K_n})=  \nonumber \\
&&\hskip 10mm \frac{\alpha}{1-\alpha } 
D_{1-\alpha}(Q_{K_1} \times \dots \times Q_{K_n} \|P_{K_1} \times \dots \times P_{K_n}),
\end{eqnarray}
and, again, for $\alpha \uparrow1$, the limits on both sides exist and coincide.
Therefore it is enough to consider either $D_\alpha(P_{K_1} \times \dots \times P_{K_n} \|Q_{K_1} \times \dots \times Q_{K_n})$ or $ D_{\alpha}(Q_{K_1} \times \dots \times Q_{K_n} \|P_{K_1} \times \dots \times P_{K_n})$.
\vskip 5mm
We first present some examples and look at special cases below. In particular,  $ D_{\pm \infty}(Q_{K_1} \times \dots \times Q_{K_n} \|P_{K_1} \times \dots \times P_{K_n})$ will be considered below.
\vskip 4mm
{\bf Examples.}
\vskip 3mm
\noindent
(i) If $K=\rho B^n_2$, then $D_\alpha(Q_K\|P_K) = D_\alpha(P_K\|Q_K) =0$ for all $- \infty \leq \alpha \leq  \infty$.
\vskip 2mm
\noindent
(ii) If $K$ is a polytope, then $\kappa_K=0$ a.e. on $\partial K$.  Thus, for $\alpha=1$,  $D_1(Q_K\|P_K) = \infty$.
For  $- \infty <\alpha <1$, $\int _{\partial K} \frac{\kappa_K^{1-\alpha} d\mu_K}{\langle x, N_K(x) \rangle ^{n-\alpha(n+1)}}=0$ and for $\alpha >1$, 
$\int _{\partial K} \frac{\kappa_K^{1-\alpha} d\mu_K}{\langle x, N_K(x) \rangle ^{n-\alpha(n+1)}}=\infty$. 
Hence 
$D_\alpha(Q_K\|P_K) = \infty$ for all $- \infty <\alpha < \infty$,  and $K$ a polytope. 
\par
Similarly, $D_1(P_K\|Q_K) = 0$ (with the convention that $0 \    \infty =0$). 
\newline
$D_\alpha(P_K\|Q_K) = - \infty$,  for  $ 1 <\alpha < \infty$ and $- \infty <\alpha <0$ and $K$ a polytope and $D_\alpha(P_K\|Q_K) = \infty$,  for  $0 <\alpha < 1$ and $K$ a polytope. 
\par
This also shows that $D_\alpha$ need not be continuous at $\alpha =1$.
\par
\noindent
For $\alpha=0$ and $\alpha = \pm \infty$, see below.
\vskip 3mm
\noindent
(iii) For $1 < r < \infty$, let 
$K=B_{r}^{n}=\{x\in\mathbb
R^{n}: \sum_{i=1}^{n}|x_{i}|^{r}\leq 1\}$ be the unit ball of $l^n_r$. We will compute $D_\alpha(Q_K\|P_K) $ and $D_\alpha(P_K\|Q_K) $ for all $- \infty < \alpha <  \infty$, $\alpha \neq 1$. The case $\alpha=1$ was considered in \cite{PW1}. The cases $\alpha=0$ and $\alpha = \pm \infty$ are treated  below.
\par
If $1 < r <2$ and $\alpha \geq \frac{1}{ 2-r}$, then $D_\alpha(P_{B_{r}^{n}} \|Q_{B_{r}^{n}}) =\infty$.
If $1 < r <2$ and $\alpha \leq -\frac{r-1}{ 2-r}$, then $D_\alpha(Q_{B_{r}^{n}} \|P_{B_{r}^{n}}) = - \infty$.
If $2 < r < \infty$ and $\alpha \leq \frac{-1}{ r-2}$, then $D_\alpha(P_{B_{r}^{n}} \|Q_{B_{r}^{n}}) = - \infty$.
If $2 < r < \infty$ and $\alpha \geq \frac{r-1}{ r-2}$, then $D_\alpha(Q_{B_{r}^{n}} \|P_{B_{r}^{n}}) =  \infty$.
In all other cases  we have
\begin{eqnarray*}
&&\hskip -9mm  D_\alpha(P_{B_{r}^{n}} \|Q_{B_{r}^{n}})=
 \frac{1}{\alpha -1}  \log \bigg[ 
  \left(\frac{\Gamma(\frac{n}{r})}{\left(\Gamma(\frac{1}{r})\right)^n}\right)^{1-\alpha}
\left(\frac{\Gamma(n(1-\frac{1}{r}))}{\left(\Gamma(1-\frac{1}{r})\right)^n}\right)^{\alpha} \\
 && \hskip 55mm \times\frac{\left(\Gamma(\frac{1-\alpha}{r} + \alpha(1-\frac{1}{r}))\right)^n}{\Gamma\left(n(\frac{1-\alpha}{r} + \alpha(1-\frac{1}{r})\right)}
\bigg]
\end{eqnarray*}
and
\begin{eqnarray*}
&&  D_\alpha(Q_{B_{r}^{n}} \|P_{B_{r}^{n}})=
 \frac{1}{\alpha -1}  \log \bigg[   \left(\frac{\Gamma(\frac{n}{r})}{\left(\Gamma(\frac{1}{r})\right)^n}\right)^{\alpha}
\left(\frac{\Gamma(n(1-\frac{1}{r}))}{\left(\Gamma(1-\frac{1}{r})\right)^n}\right)^{1-\alpha} \\
 && \hskip 55mm \times
\  \frac{\left(\Gamma(\frac{\alpha}{r} + (1-\alpha)(1-\frac{1}{r}))\right)^n}{\Gamma\left(n(\frac{\alpha}{r} + (1-\alpha)(1-\frac{1}{r})\right)}.
\bigg]
\end{eqnarray*}
\vskip 5mm
Now we introduce $L_p$-affine surface areas for a convex body $K$ in $\mathbb{R}^n$.
$L_p$-affine surface area, an extension of affine surface area, 
was introduced by Lutwak in the ground
breaking paper \cite{Lu2}  for $p >1$ and for general $p$ by Sch\"utt and Werner \cite{SW2004}.
For real  $p \neq -n$, we define  the
$L_p$-affine surface area $as_{p}(K)$ of $K$ as in \cite{Lu2} ($p
>1$) and \cite{SW2004} ($p <1, p \neq -n$) by
\begin{equation} \label{def:paffine}
as_{p}(K)=\int_{\partial K}\frac{\kappa_K(x)^{\frac{p}{n+p}}}
{\langle x,N_{ K}(x)\rangle ^{\frac{n(p-1)}{n+p}}} d\mu_{ K}(x) 
\end{equation}
and
\begin{equation}\label{def:infty}
as_{\pm\infty}(K)=\int_{\partial K}\frac{\kappa _K (x)}{\langle
x,N_{K} (x)\rangle ^{n}} d\mu_{K}(x), 
\end{equation}
provided the above integrals exist.
In particular, for $p=0$
$$
as_{0}(K)=\int_{\partial K} \langle x,N_{ K}(x)\rangle
\,d\mu_{K}(x) = n|K|.
$$
The case $p=1$ is the classical affine surface area which goes
back to Blaschke. It is  independent
of the position of $K$ in space.
$$
as_{1}(K)=\int_{\partial K} \kappa_{ K}(x)^\frac{1}{n+1} 
\,d\mu_{K}(x).
$$
Originally a basic affine invariant
from the field of affine differential geometry, it has recently
attracted increased attention too (e.g. 
\cite{LR1,  Lu2, MW1, SW1,  W1994}).
\vskip 3mm
If $K$ is in $C^2_+$, then $d\mu_K= f_K d \sigma$ and  then the $L_p$-affine surface areas, for all  $p\neq -n$, can  be written as
\begin{equation}\label{sphereform}
as_p(K)=\int_{S^{n-1}}\frac{f_{K}(u)^{\frac{n}{n+p}}}
{h_K(u)^{\frac{n(p-1)}{n+p}}} d\sigma(u).
\end{equation}
In particular,
$$
as_{\pm\infty}(K)=\int_{S^{n-1}}\frac{d\sigma(u)}
{h_K(u)^{n}}  = n |K^\circ|.
$$
Recall that $f_{K}(u)$ is the curvature function of
$K$ at
 $u$, i.e., the reciprocal of the Gauss
curvature $\kappa _K(x)$ at this point $x \in
\partial K$, the boundary of $K$,  that has $u$ as its outer normal. 
\vskip 3mm
The mixed
$p$-affine surface area, $as_{p}(K_1, \cdots,K_{n})$, of $n$
convex bodies  $K_i \in C^2_+$ was introduced - for $p\geq 1$  in  \cite{Lu1} and extended to all $p$ in \cite{WY2010} - as
\begin{equation}\label{def:pmixed}
as_p(K_1, \cdots, K_{n})=\int
_{S^{n-1}}\bigg[h_{K_1}(u)^{1-p}f_{K_1}(u)\cdots
h_{K_n}^{1-p}f_{K_n}(u)\bigg]^{\frac{1}{n+p}}\,d\sigma (u).
\end{equation}

\vskip 5mm
Then we observe the following remarkable fact which connects  $L_p$-Brunn Minkowki theory and information theory: 
\newline
$L_p$-affine surface areas of a convex body  are Hellinger integrals - or   exponentials of R\'enyi divergences - of the cone measures of $K$ and $K^\circ$. For $\alpha=1$, such a connection was already observed in \cite{PW1}, namely
\begin{equation}\label{omega-lim}
\frac{|K|}{|K^\circ|}e^{-  D_{KL}(P_K\|Q_K)} =  \lim_{p \rightarrow \infty} \left(\frac{as_p(K)}{n |K^\circ|}\right)^\frac{n+p}{n}.
\end{equation}
\vskip 4mm
Now we have more generally
\begin{theorem}\label{p-aff=renyi}
Let $K$ be a convex body in $\mathbb {R}^n$. Let $- \infty < \alpha < \infty$. $\alpha \neq1$.
 Then
\begin{equation*}\label{p-aff=renyi1}
D_\alpha(P_K\|Q_K)=  \frac{1}{\alpha -1} \log\left( \frac{as_{n\frac{\alpha}{1-\alpha}}(K)}{n |K|^{1-\alpha} |K^\circ|^{\alpha} }\right).
\end{equation*}
\begin{equation*}
D_\alpha(Q_K\|P_K)=  \frac{1}{\alpha -1} \log\left( \frac{as_{n\frac{1-\alpha}{\alpha}}(K)}{n |K|^\alpha |K^\circ|^{1-\alpha} }\right).
\end{equation*}
Equivalently, for  all $- \infty \leq p \leq \infty$, $p \neq -n$, 
\begin{eqnarray*}\label{p-aff=renyi2}
\frac{as_p(K)}{n |K|^\frac{n}{n+p}   |K^\circ| ^\frac{p}{n+p}} &=&\text{Exp}\left(-\frac {n}{n+p}  D_{\frac{p}{n+p}}\left(P_K\|Q_K\right) \right)\\
&=&\text{Exp}\left( -\frac {p}{n+p}  D_{\frac{n}{n+p}}(Q_K\|P_K)\right)
\end{eqnarray*}
In particular,
\begin{eqnarray*}\label{p-aff=renyi2}
\frac{as_1(K)}{ n |K|^\frac{n}{n+1}   |K^\circ| ^\frac{1}{n+1}} &=& \text{Exp}\left( -\frac {n}{n+1}  D_{\frac{1}{n+1}}(P_K\|Q_K) \right)\\
&=&\text{Exp}\left(  -\frac {1}{n+1}  D_{\frac{n}{n+1}}(Q_K\|P_K) \right).
\end{eqnarray*}

\end{theorem}
\vskip 4mm

{\bf Remarks.}
 \vskip 2mm
 \noindent 
(i) Theorem \ref{p-aff=renyi}  can also be written as
\begin{eqnarray*}
\left(\frac{as_p(K)}{n   |K^\circ| }\right)^\frac{n+p}{n}= \frac{ |K|}  { |K^\circ| }\ e^{-  D_{\frac{p}{n+p}}(P_K\|Q_K)}.
\end{eqnarray*}
If we  now let $p \rightarrow \infty$, we recover (\ref{omega-lim}).
Also from Theorem \ref{p-aff=renyi} 
\begin{eqnarray*}
\left(\frac{as_p(K)}{n   |K| }\right)^\frac{n+p}{p}= \frac{ |K^\circ|}  { |K| }\  e^{-  D_{\frac{n}{n+p}}(Q_K\|P_K)}.
\end{eqnarray*}
If we  let $p \rightarrow 0$, then we get 
\begin{equation}\label{Alpha}
\lim_{p \rightarrow 0} \left(\frac{as_p(K)}{n |K|}\right)^\frac{n+p}{p} = \frac{ |K^\circ|}  { |K| } \ e^{-  D_{KL}(Q_K\|P_K)}.
\end{equation}
We will comment on  these expressions  in Section 3.
 \vskip 2mm
 \noindent 
 (ii) 
If $- \infty  <  \alpha \leq 0$, then $- \infty \leq   p=n\frac{1-\alpha}{\alpha} < -n$. Thus, for this range of $\alpha$,  we get the $L_p$-affine surface area in the range smaller than $-n$.
If $0 \leq \alpha < \infty$, then $-n  < p=n\frac{1-\alpha}{\alpha} \leq \infty$. Thus, for this range of $\alpha$, we get the $L_p$-affine surface area in the range greater than $-n$. In particular, for $0 \leq \alpha \leq 1$, we get the $L_p$-affine surface area for $0 \leq p \leq \infty$. 
\par
If $- \infty \leq   \alpha< 1$, then $-n  <  p=n\frac{\alpha}{1-\alpha} \leq \infty$. Thus, for this range of $\alpha$, we get the $L_p$-affine surface area in the range greater $-n$.
If $1 < \alpha \leq  \infty$, then $- \infty \leq  p=n\frac{\alpha}{1-\alpha} < -n$. Thus, for this range of $\alpha$, we get the $L_p$-affine surface area in the range smaller than $-n$.

\vskip 5mm
\noindent
\begin{theorem}\label{mixed-aff=renyi}
Let $K_1, \dots K_n$ be convex bodies in $C^2_+$. Then, for all $\alpha \neq 1$
\begin{equation*}
D_\alpha(P_{K_1} \times \dots \times P_{K_n} \|Q_{K_1} \times \dots \times Q_{K_n})=  \frac{1}{\alpha -1} \log\left( \frac{as_{n\frac{\alpha}{1-\alpha}}(K_1, \dots, K_n)}{n \prod_{i=1}^n |K_i|^\frac{1-\alpha}{n} |K_i^\circ|^\frac{\alpha}{n} }\right)
\end{equation*}
and
\begin{equation*}
D_\alpha(Q_{K_1} \times \dots \times Q_{K_n} \|P_{K_1} \times \dots \times P_{K_n})=  \frac{1}{\alpha -1} \log\left( \frac{as_{n\frac{1-\alpha}{\alpha}}(K_1, \dots, K_n)}{n \prod_{i=1}^n |K_i|^\frac{\alpha}{n} |K_i^\circ|^\frac{1-\alpha}{n} }\right).
\end{equation*}

\end{theorem}

\vskip 4mm
{\bf Remark.}
 \vskip 2mm
 \noindent 
The expressions in Theorem \ref{mixed-aff=renyi}  can also be written as
\begin{eqnarray*}
\left(\frac{as_{n\frac{\alpha}{1-\alpha}}(K_1, \dots, K_n)}{n \prod_{i=1}^n  |K_i^\circ| ^\frac{1}{n}}\right)^\frac{1}{1-\alpha}= \prod_{i=1}^n \left( \frac{ |K_i|}  { |K_i^\circ | } \right)^\frac{1}{n} \ e^{-  D_\alpha(P_{K_1} \times \dots \times P_{K_n} \|Q_{K_1} \times \dots \times Q_{K_n})}.
\end{eqnarray*}
and
\begin{eqnarray*}
\left(\frac{as_{n\frac{1-\alpha}{\alpha}}(K_1, \dots, K_n)}{n \prod_{i=1}^n  |K_i |^\frac{1}{n}}\right)^\frac{1}{1-\alpha}= \prod_{i=1}^n \left( \frac{ |K_i ^\circ |}  { |K_i | } \right)^\frac{1}{n} \ e^{-  D_\alpha(Q_{K_1} \times \dots \times Q_{K_n} \|P_{K_1} \times \dots \times P_{K_n})}.
\end{eqnarray*}
\par
If we  now let in the first expression  $\alpha \rightarrow 1$ respectively, putting $p=n \frac{\alpha}{1-\alpha}$,  $p \rightarrow \infty$, we  get
\begin{eqnarray}\label{mixed1=renyi,1}
&&\hskip -10mm \prod_{i=1}^n \left( \frac{ |K_i|}  { |K_i^\circ| } \right)^\frac{1}{n} \ e^{-  D_1(P_{K_1} \times \dots \times P_{K_n} \|Q_{K_1} \times \dots \times Q_{K_n})}  \nonumber \\
&& \hskip 10mm = \lim_{\alpha \rightarrow 1}  \left(\frac{as_{n\frac{\alpha}{1-\alpha}}(K_1, \dots, K_n)}{n \prod_{i=1}^n  |K_i^\circ| ^\frac{1}{n}}\right)^\frac{1}{1-\alpha} \nonumber \\
&& \hskip 10mm= \lim _{p \rightarrow \infty }  \left(\frac{as_{p}(K_1, \dots, K_n)}{n \prod_{i=1}^n  |K_i^\circ| ^\frac{1}{n}}\right)^\frac{n+p}{n}.
\end{eqnarray}
\par
If we   let in the second expression $\alpha \rightarrow 1$, respectively, putting $p=n \frac{1-\alpha}{\alpha}$,  $p \rightarrow 0$, we  get
\begin{eqnarray}\label{mixed1=renyi,2}
&&\hskip -10mm \prod_{i=1}^n \left( \frac{ |K_i ^\circ|}  { |K_i } \right)^\frac{1}{n} \ e^{-  D_1(Q_{K_1} \times \dots \times Q_{K_n} \|P_{K_1} \times \dots \times P_{K_n})}  \nonumber \\
&& \hskip 10mm = \lim_{\alpha \rightarrow 1}  \left(\frac{as_{n\frac{1-\alpha}{\alpha}}(K_1, \dots, K_n)}{n \prod_{i=1}^n  |K_i| ^\frac{1}{n}}\right)^\frac{1}{1-\alpha} \nonumber \\
&& \hskip 10mm= \lim _{p \rightarrow 0}  \left(\frac{as_{p}(K_1, \dots, K_n)}{n \prod_{i=1}^n  |K_i| ^\frac{1}{n}}\right)^\frac{n+p}{p}.
\end{eqnarray}
 \vskip 2mm
 \noindent 
We will comment on  these quantities in Section 3.

\vskip 4mm
{\bf Special Cases.}
\vskip 2mm
\noindent
(i) If $\alpha = \frac{1}{2}$, then 
\begin{equation*}\label{a=1/2,1}
  D_{\frac{1}{2}}(Q_K\|P_K) =   D_{\frac{1}{2}}(P_K\|Q_K) =  - 2 \log\left( \frac{ as_{n}(K)}{  n |K|^\frac{1}{2} |K^\circ| ^\frac{1}{2}} \right),
  \end{equation*}
 and
 $\frac{ as_{n}(K)}{  n |K|^\frac{1}{2} |K^\circ| ^\frac{1}{2}}$ is the Bhattcharyya coefficient of $p_K$ and $q_K$.
 
\begin{eqnarray*}\label{a=/2,2}
 &&\hskip -15mm  D_{\frac{1}{2}}(Q_{K_1} \times \dots \times Q_{K_n} \|P_{K_1} \times \dots \times P_{K_n})= \\  
 && \hskip 20mm  = D_{\frac{1}{2}}(  P_{K_1} \times \dots \times P_{K_n} \|Q_{K_1} \times \dots \times Q_{K_n})\\
 && \hskip 35mm =   - 2 \log\left( \frac{ as_{n}(K_1, \dots, K_n)}{n \prod_{i=1}^n  |K_i|^\frac{1}{2n} |K_i^\circ| ^\frac{1}{2n}} \right)
  \end{eqnarray*}
\vskip 2mm
\noindent
(ii)
If $\alpha = 0$, then $D_{0}(P_K\|Q_K)= 0$. Likewise, 
\begin{equation}\label{a=0}
  D_{0}(Q_K\|P_K)=   - \log\left( \frac{ as_{\infty}(K)}{  n |K^\circ|} \right)
  \end{equation}
  which, if $K$ is sufficiently smooth, is equal to 
  $$ 
 - \log\left( \frac{ as_{\infty}(K)}{  n |K^\circ|} \right)
 = - \log\left( \frac{\int_{\partial K} \frac{\kappa_K(x) d \mu(x)}{\langle x , N_K(x) \rangle^n}}{n |K^\circ|} \right)=  -\log 1 =0
  $$
  and equal to $\infty$ if $K$ is a polytope.
\begin{eqnarray*}\label{a=/2,2}
D_{0}(  P_{K_1} \times \dots \times P_{K_n} \|Q_{K_1} \times \dots \times Q_{K_n}) = - \log\left( \frac{ as_{0}(K_1, \dots, K_n)}{n \prod_{i=1}^n  |K_i|^\frac{1}{n} } \right)
  \end{eqnarray*}
and
\begin{eqnarray*}\label{a=/2,2}
  D_{0}(Q_{K_1} \times \dots \times Q_{K_n} \|P_{K_1} \times \dots \times P_{K_n})&=&  
   -  \log\left( \frac{ as_{\infty}(K_1, \dots, K_n)}{n \prod_{i=1}^n  |K_i^\circ| ^\frac{1}{n}} \right) \\
 &=&  -  \log\left(\frac{ \tilde{V}(K_1, \dots, K_n)}{ \prod_{i=1}^n  |K_i^\circ| ^\frac{1}{n}} \right),
  \end{eqnarray*}
where $\tilde{V}(K_1, \dots, K_n)$ is the dual mixed volume introduced by Lutwak in \cite{Lut1975}.
\vskip 2mm
\noindent
(iii) 
If $ \alpha \rightarrow  \infty$, then $p=n\frac{1-\alpha}{\alpha} \rightarrow -n$ from the right. Therefore,  by definition,
$D_{\infty}(Q_K\|P_K) =  \log   \left(\sup_x \text{ess} \  \frac{q_K(x)}{p_K(x)}\right) =  \log  \left(\sup_x \text{ess} \ \frac{\langle x, N_{K}(x) \rangle^{n+1}|K^\circ | } { \kappa_{K}(x)  |K|} \right)$. On the other hand
\begin{eqnarray*}
\lim_{\alpha \rightarrow \infty}\left( \frac{as_{n\frac{1-\alpha}{\alpha}}(K)}{n |K| ^\alpha |K^\circ|^{1-\alpha}}\right)^\frac{1}{\alpha -1} &=  &
\frac{|K^\circ|}{|K|} \lim _{\alpha \rightarrow  \infty} \bigg\|  \frac{\langle x , N_K(x) \rangle^{n+ \frac{\alpha}{1-\alpha}} }{\kappa_K(x) } \bigg\|_{L_{\alpha-1}} \\
&=& \frac{|K^\circ|}{|K|} \bigg\| \frac{\langle x , N_K(x) \rangle^{n+1}}{\kappa_K(x)} \bigg\|_{L_{\infty}},
\end{eqnarray*}
which is thus 
consistent with the definition of $D_{\infty}(Q_K\|P_K)$. 
Similarly, one shows that,  if $ \alpha \rightarrow  \infty$, then $p=n\frac{\alpha}{1-\alpha} \rightarrow -n$ from the left.
Hence,  by definition, $D_{\infty}(P_K\|Q_K)= \log  \left(\sup_x \text{ess} \  \frac{q_K(x)}{p_K(x)}\right) =  \log   \left(\sup_x \text{ess} \  \frac{\kappa_{K}(x)  |K|}{\langle x, N_{K}(x) \rangle^{n+1}|K^\circ | }\right)$, which is consistent with $\lim_{\alpha \rightarrow \infty}\left( \frac{as_{n\frac{\alpha}{1-\alpha}}(K)}{n |K| ^{1-\alpha} |K^\circ|^{\alpha}}\right)^\frac{1}{\alpha -1}$.
\vskip 4mm
\noindent
Thus, also it would make most sense to define
\begin{equation}\label{def:-n+}
\lim_{p \rightarrow -n^+} as_{p}(K) = \sup_{x \in \partial K} \text{ess} \  \frac{\langle x , N_K(x) \rangle^{n+1}}{\kappa_K(x)}.
\end{equation}
and
\begin{equation}\label{def:-n-}
\lim_{p \rightarrow -n^-} as_{p}(K) =  \sup_{x \in \partial K} \text{ess} \   \frac{\kappa_K(x)}{\langle x , N_K(x) \rangle^{n+1}},
\end{equation}
which would imply that $\lim_{p \rightarrow -n} as_{p}(K)$ does not exist.
\vskip 3mm
\noindent
If $ \alpha \rightarrow  -\infty$, then $p=n\frac{1-\alpha}{\alpha} \rightarrow -n$ from the left and by (\ref{minusinfty}), 
$D_{-\infty}(Q_K\|P_K) = - D_{\infty}(P_K\|Q_K) $. On the other hand, 
\begin{eqnarray*}
\lim_{\alpha \rightarrow \infty}\log\left( \frac{as_{n\frac{1-\alpha}{\alpha}}(K)}{n |K| ^\alpha |K^\circ|^{1-\alpha}}\right)^\frac{1}{\alpha -1} &= &
\log \left(\frac{1}{\sup_{x} \frac{\kappa_{K}(x)  |K|}{\langle x, N_{K}(x) \rangle^{n+1}|K^\circ | }}\right)\\
&=& - \log \left(\sup_{x} \frac{\kappa_{K}(x)  |K|}{\langle x, N_{K}(x) \rangle^{n+1}|K^\circ | }\right) \\
&=& - D_{\infty}(P_K\|Q_K), 
\end{eqnarray*} 
hence this is also consistent with the definitions.
Similar considerations  hold for  $D_{-\infty}(P_K\|Q_K)$ and  $D_{\alpha}(  P_{K_1} \times \dots \times P_{K_n} \|Q_{K_1} \times \dots \times Q_{K_n})$ and $D_{\alpha}(  Q_{K_1} \times \dots \times Q_{K_n} \|P_{K_1} \times \dots \times P_{K_n})$.

\vskip 5mm

Having identified $L_p$-affine surface areas as R\'enyi divergences, 
we can now translate known results from one theory to the other.
\par
Affine invariance of $L_p$-affine surface areas translates into affine invariance of R\'enyi divergences:  
For all $p \neq -n$, 
$as_p(T(K)) = | \text{det} \ T|^\frac{n-p}{n+p} as_p(K)$ (see \cite{SW2004}). Theorem  \ref{p-aff=renyi}  then implies that
for all linear maps $T$ with $\text{det} \ T \neq 0$, for all $- \infty < \alpha < \infty$, $\alpha \neq 1$,
\begin{eqnarray*}
D_\alpha (P_{T(K)} \| Q_{T(K)}) =  D_{\alpha} (P_{K} \| Q_{K}) 
\end{eqnarray*}
and
\begin{eqnarray*}
D_\alpha (Q_{T(K)} \| P_{T(K)}) =  D_{\alpha} (Q_{K} \| P_{K}).
\end{eqnarray*}

The case $\alpha=1$  was treated in \cite{PW1}.
\vskip 2mm
\noindent
As $as_p(T(K_1), \dots, T(K_n) ) = | \text{det} \ T|^\frac{n-p}{n+p} as_p(K_1, \dots, K_n)$ (see \cite{WY2010}), it follows from Theorem  \ref{mixed-aff=renyi} that for all linear maps $T$ with $\text{det} \ T \neq 0$, for all $- \infty < \alpha < \infty$, $\alpha \neq 1$,
\begin{eqnarray*}
&&\hskip -15mm D_{\alpha}(  P_{T(K_1)} \times \dots \times P_{T(K_n)} \|Q_{T(K_1)} \times \dots \times Q_{T(K_n)})\\
&&\hskip 15mm =  D_{\alpha}(  P_{K_1} \times \dots \times P_{K_n} \|Q_{K_1} \times \dots \times Q_{K_n}) 
\end{eqnarray*}
and
\begin{eqnarray*}
&&\hskip -15mm D_{\alpha}(  Q_{T(K_1)} \times \dots \times Q_{T(K_n)} \|P_{T(K_1)} \times \dots \times P_{T(K_n)})\\
&& \hskip 15mm = D_{\alpha}(  Q_{K_1} \times \dots \times Q_{K_n} \|P_{K_1} \times \dots \times P_{K_n}).
\end{eqnarray*}
The case $\alpha=1$  is  in \cite{PW1}.

\vskip 3mm

Moreover, all inequalities and results mentioned in e.g. \cite{WY2008} about $L_p$-affine surface area and  in e.g. \cite{WY2010} about mixed $L_p$-affine surface area can be translated into the corresponding inequalities and results about R\'enyi divergences. Conversely, results about R\'enyi divergences from e.g. \cite{Erven/Harre} have consequences for $L_p$-affine surface areas. We mention only a few.
\vskip 4mm
\begin{proposition} Let $K$ be a convex body in $C^2_+$.
\vskip 2mm
\noindent
(i) Then for all $-\infty  \leq  \alpha \leq  \infty$,
$$
(1-\alpha) D_\alpha (Q_{K^\circ} \| P_{K^\circ}) = \alpha D_{1-\alpha} (Q_{K} \| P_{K})
$$
and
$$
(1-\alpha) D_\alpha (P_{K^\circ} \| Q_{K^\circ}) = \alpha D_{1-\alpha} (P_{K} \| Q_{K}) 
$$
The equalities hold trivially if $\alpha=0$ or $\alpha=1$.
\vskip 2mm
\noindent
(ii) Let $K_i$, $1 \leq i \leq n$,  be  convex bodies in $C^2_+$. Then for all $ 0 \leq \alpha$
\begin{eqnarray*}
as_{n\frac{\alpha}{1-\alpha}}(K_1, \dots, ,K_n) &=&\int_{S^{n-1}} \prod_{i=1}^m \left[f_{K_i} h_{K_i}^{1-\frac{n \alpha}{1-\alpha}} \right] ^\frac{1-\alpha}{n} d \sigma \\ &=&  
\prod_{i=1}^m
\int_{S^{n-1}}\left[f_{K_i} h_{K_i}^{1-\frac{n\alpha}{1-\alpha} }\right] ^\frac{1-\alpha}{n} d \sigma, 
\end{eqnarray*}
i.e. we can interchange integration and product.
\vskip 2mm
\noindent
(iii) 
Let $K$ and $L$  be  convex bodies in $C^2_+$. Let $0 \leq p \leq \infty$. Let $0 \leq \lambda \leq 1$.
Then
\begin{eqnarray*}
\int_{S^{n-1}} \left[\lambda \frac{f_{K} h_{K}}{|K| } + (1-\lambda) \frac{f_{L} h_{L}}{|L| } \right] ^\frac{n}{n+p} 
\left[ \frac{\lambda} {h_{K}^n |K^\circ| } + \frac{1-\lambda}{ h_{L}^n |L^\circ| } \right] ^\frac{p}{n+p}d \sigma\\
\geq \left( \frac{as_p(K)}{|K|^\frac{n}{n+p} |K^\circ |^\frac{p}{n+p}}\right)^\lambda \  \left( \frac{as_p(L)}{|L|^\frac{n}{n+p} |L^\circ |^\frac{p}{n+p}}\right)^{1-\lambda}
\end{eqnarray*}
with equality iff $K=L$. Equality holds trivially if $p=0$ or $p= \infty$ or $\lambda=0$ or $\lambda=1$.

\end{proposition}
\vskip 2mm
\noindent
{\bf Proof.}
\par
\noindent
(i) For $-\infty < \alpha < \infty$, (i) follows from the duality formula 
$as_p(K)= as_{\frac{n^2}{p}}(K^\circ)$, or, formulated in a more symmetric way, using the parameter $\alpha=\frac{p}{n+p}$
$$
as_{n\frac{\alpha}{1-\alpha}} (K)= as_{n\frac{1-\alpha}{\alpha}}(K^\circ).
$$
This identity was  proved for $p>0$ in \cite{Hug} and - with a different proof - for all other $p$ in \cite{WY2008}. 
\par
Let now  $\alpha = \infty$. Then,  on the one hand
\begin{equation}\label{gl1}
\lim_{\alpha \rightarrow \infty} \frac{1-\alpha}{\alpha}  D_\alpha (Q_{K^\circ} \| P_{K^\circ}) = - D_\infty (Q_{K^\circ} \| P_{K^\circ}) = - \log \sup_{x \in \partial K^\circ} \text{ess} \  \frac{q_{K^\circ}(x)}{p_{K^\circ}(x)}.
\end{equation} 
On the other hand, by (\ref{K-renyi5}), 
\begin{equation}\label{gl2}
D_{-\infty}  (Q_{K} \| P_{K})= - D_{-\infty}  (P_{K} \| Q_{K}) = -\log \sup_{x \in \partial K} \text{ess} \ \frac{p_{K}(x)}{q_{K}(x)}.
\end{equation}
 (\ref{gl1}) equals (\ref{gl2}), as (see \cite{Hug}) for $x \in \partial K$, $y \in \partial K^\circ $ such that 
$\langle x,y \rangle =1$,  
\begin{equation*} 
  \langle y, N_{K^\circ} (y)  \rangle \langle  x,
N_{K} (x) \rangle =\left(\kappa _{K^\circ}(y)\kappa
_{K}(x)\right)^{\frac{1}{n+1}}.
\end{equation*}
Similarly, for $\alpha =- \infty$.
\vskip 2mm
\noindent
(ii) follows from Theorem \ref{mixed-aff=renyi} and the fact that \cite{Erven/Harre}
$$
D_{\alpha} (Q_{K_1} \times \dots \times Q_{K_n} \| P_{K_1} \times \dots \times P_{K_n})
= \sum_{i=1}^n D_{\alpha} (Q_{K_i}  \| P_{K_i}), 
$$ respectively 
the corresponding equation for $D_{\alpha} (P_{K_1} \times \dots \times P_{K_n} \| Q_{K_1} \times \dots \times Q_{K_n}) $.
\vskip 2mm
\noindent
(iii) 
For $0 \leq \alpha \leq 1$, $D_{\alpha}  (Q_{K} \| P_{K})$,  respectively $D_{\alpha}  (P_{K} \| Q_{K})$,  are jointly convex \cite {Erven/Harre}.  We put $p=n \frac{1-\alpha}{\alpha}$ respectively $p=n \frac{\alpha}{1-\alpha}$ and use the joint convexity together with Theorem \ref {p-aff=renyi}.
\par
If $p \neq 0, \infty$ and $\lambda \neq 0,1$, then equality implies that $K=L$ as the logarithm is strictly concave.

\vskip 5mm
\section{Geometric interpretation of R\'enyi Divergence}

In this section we present geometric interpretations of 
R\'enyi divergences  $D_\alpha$   
of convex bodies, for all $\alpha$.  Geometric interpretations for the case $\alpha=1$, the relative entropy,   were given first
 in \cite{PW1}  in terms of $L_p$-centroid bodies. 
Recall that  for a convex body  $K$  in $\mathbb R^n$ of volume $1$ and $1\leq p \leq \infty$,  the $L_{p}$-centroid body $Z_{p}(K)$ is  this convex body that has support function 
\begin{equation*} \label{def:Zp}
h_{Z_{p}(K)}(\theta) = \left(\int_{K}|\langle x, \theta\rangle |^{p} dx \right)^{1/p}.
\end{equation*}
\par
 Now that we observed that R\'enyi divergences are logarithms of 
$L_p$-affine surface areas, we can use their geometric characterizations to obtain the ones for  R\'enyi divergences. 
We will mostly concentrate on the geometric characterization of $L_p$-affine surface areas via the {\em surface bodies} \cite{SW2004}
and {\em illumination surface bodies} \cite{WY2010},  though there are many more available (see e.g. \cite{MW2,
SW2002, W2, WY2008})
\par
Even more is  gained. Firstly, we need not assume that the body is symmetric as in \cite{PW1} nor that it has $C^2_+$ boundary as it was needed in \cite{PW1}, to obtain the desired 
geometric interpretation for the $D_\alpha$ for all $\alpha$. Weaker regularity assumptions on the boundary suffice.
\par
Secondly, in the context of the  $L_p$-centroid bodies, the relative entropies appeared only after performing  a second order expansion of certain expressions. Now, using the surface bodies or illumination surface bodies, already a first order expansion makes them appear. Thus,  these bodies detect ``faster"  details  of the  boundary of a convex body than the  $L_p$-centroid bodies.

\vskip 5mm
Let $K$ be a convex body in $\mathbb R^{n}$.
Let $f:\partial K\rightarrow \Bbb R$  be a nonnegative,  integrable, function. 
Let $s\geq 0$. 
\par
The {\em surface body}  $K_{f,s}$,  introduced in \cite{SW2004},  is the intersection of all closed half-spaces
$H^{+}$
whose defining hyperplanes $H$ cut off a set of 
$f \mu_K$-measure less
than or equal to $s$ from $\partial K$. More precisely,
\begin{equation*}\label{def1}
K_{f,s}= \bigcap _{\int_{\partial K \cap H^-}f d\mu_K \leq s} H^+.
\end{equation*}
\vskip 3mm
\noindent
The {\em illumination surface body}  $K^{f,s}$ \cite {WY2010}  is defined as
$$\displaystyle K^{f,s}=\left\{x: \mu _f({\partial K \cap
\overline{[x,K]\backslash K}})\leq s\right\},$$ 
where for sets $A$ and $B$ (respectively points $x$ and $y$) in $\mathbb{R}^n$,   $[A,B]=
\{\lambda a + 1-\lambda b: a \in A, b \in B, 0 \leq \lambda \leq 1\}$ (respectively  $[x,y] =
\lambda x + 1-\lambda y:  0 \leq \lambda \leq 1\}$) is the convex hull of $A$ and $B$ (respectively $x$ and $y$).

\vskip 5mm
For $ x \in \partial K$ and $s > 0$ and $f$ and $K_{f,s}$ as above, we put
$$x_s=[0,x] \cap \partial K_{f,s}.
$$

The {\em minimal function} $M_f : \partial K \rightarrow \mathbb{R}$
\begin{equation}\label{min}
M_{f}(x)=\inf_{0<s} \  \frac{\int_{\partial K\cap
H^{-}(x_{s},N_{ K_{f,s}}(x_{s}))}f\ d\mu_{ K}}{\mbox{vol}_{n-1}\left(\partial K\cap
H^{-}(x_{s},N_{ K_{f,s}}(x_{s}))\right)}
\end{equation}
was introduced  in   \cite{SW2004}.  
$H(x,\xi)$ is the hyperplane through $x$ and orthogonal to
$\xi$.
$H^{-}(x,\xi)$ is the closed halfspace containing the point $x+\xi$,
$H^{+}(x,\xi)$ the other halfspace.

For $x \in \partial K$, we define $r(x)$ as
the maximum of all real numbers $\rho$ so that
$B_{2}^{n}(x-\rho N_{ K}(x),\rho)\subseteq K$. 
Then we formulate an  integrability condition for the minimal function
\begin{equation}\label{intcond}
\int_{\partial K}\frac{d\mu_{
K}(x)}{\left((M_{f}(x)\right)^{\frac{2}{n-1}}r(x)}
<\infty.
\end{equation}

\vskip 4mm
The following theorem was proved in \cite{SW2004}.
\vskip 3mm
\begin{theorem}\label{thm:surface}  \cite{SW2004}
Let $K$ be a convex body in $\mathbb R^{n}$.
Suppose that $f:\partial K\rightarrow \mathbb R$ is an integrable,
almost everywhere strictly positive function
that satisfies the integrability condition (\ref{intcond}).
Then
$$
c_{n}\lim_{s \to 0}
\frac{|K|-|K_{f,s}|}
{s^\frac{2}{n-1}}=
\int_{\partial K} \frac{\kappa^\frac{1}{n-1}}
{f^\frac{2}{n-1}}d\mu_{K}.
$$
 $c_n=2 |B_2^{n-1}|^{\frac{2}{n-1}}$. 
\end{theorem}

\vskip 4mm
Theorem \ref{thm:surface}  was used  in \cite{SW2004}  to give  geometric interpretations of $L_p$-affine surface area.
Now we use this theorem to give 
geometric interpretations of R\'enyi divergence of order $\alpha$ for all $\alpha$ for cone measures of convex bodies.
First  we treat the case $\alpha \neq 1$.
\vskip 4mm
\begin{corollary}\label{cor:asp}
Let $K$ be a convex body in $\mathbb R^{n}$.
\vskip 2mm
\noindent
 For $- \infty \leq p \leq \infty$, $p \neq -n$,  let $f_p: \partial K \rightarrow \mathbb{R}$  be defined as
\begin{equation*}
f_p(x)= \frac{\langle x,N_{ K}(x)\rangle ^\frac{(n-1)n(p-1)}{2(n+p)}}{\kappa_K(x)^\frac{n(p-1)-2p}{2(n+p)}}.
\end{equation*}
 If  $f_p$ is almost everywhere strictly positive and satisfies the integrability condition (\ref {intcond}), then
 \begin{eqnarray*}
\frac{c_n }{n |K|^{\frac{n}{n+p}}   |K^\circ| ^\frac{p}{n+p} }\   \lim_{s \to 0}
\frac{|K|
-|K_{f_p,s}|}{  s^\frac{2}{n-1}}
=  \text{Exp} \left( - \frac {p}{n+p} D_{\frac{n}{n+p}}(Q_K\|P_K) \right), 
\end{eqnarray*}
\vskip 2mm
\noindent
and , provided $p \neq \pm \infty$, 
\begin{eqnarray*}
\frac{c_n }{n |K|^{\frac{n}{n+p}}   |K^\circ| ^\frac{p}{n+p} }\   \lim_{s \to 0}
\frac{|K|
-|K_{f_p,s}|}{  s^\frac{2}{n-1}}
=  \text{Exp} \left( - \frac {n}{n+p}  D_{\frac{p}{n+p}}(P_K\|Q_K)\right).
\end{eqnarray*}
\vskip 2mm
\noindent
If $K$ is in $C^2_+$,  the last equation also holds for $p = \pm \infty$.
\end{corollary}
\vskip 4mm
\noindent
{\bf Proof.} 
The proof of the corollary follows immediately from Theorems \ref{thm:surface} and \ref{p-aff=renyi}.
\vskip 5mm
The next corollary treats the case $\alpha =1$. There, we need to make additional regularity assumptions on the 
boundary of $K$. Those are weaker though than $C^2_+$. 
\vskip 5mm
\begin{corollary}\label{cor:dkl}
Let $K$ be a convex body in $\mathbb R^{n}$.  Assume that $K$ is such that there are $0 < r \leq R < \infty$ so that for all $x \in \partial K$
\begin{equation}\label{condition}
B^n_2(x-r N_K(x), r) \subset K \subset  B^n_2(x-R N_K(x), R).
\end{equation}
Let $f_{PQ}: \partial K \rightarrow \mathbb{R}$ and $f_{QP}: \partial K \rightarrow \mathbb{R}$ be defined by
$$
f_{PQ}(x) = 
\frac {\left(n |K^\circ|  \langle x, N_{K}(x) \rangle \right) ^{\frac{n-1}{2}}} {\kappa_K(x)^\frac{n-2}{2}}
\left( \log \left(
\frac{ R^{2n}|K| \ \kappa_K(x)} {r^{2n} |K^\circ| \  \langle x, N_{K}(x) \rangle^{n+1}} \right) \right)^{-\frac{n-1}{2}}, 
$$
$$
f_{QP}(x) = \left(
\frac {n |K| }{ \langle x, N_{K}(x) \rangle }\right)^\frac{n-1}{2}
\kappa_K(x)^{\frac{1}{2}}  
\left( \log \left(
\frac{R^{2n}|K^\circ | \  \langle x, N_K(x) \rangle ^{n+1}} {r^{2n} |K| \ \kappa_K(x)} \right) \right)^{-\frac{n-1}{2}}.
$$

Then $f_{PQ}$  and  $f_{QP}$ are almost everywhere strictly positive, satisfy the integrability condition (\ref {intcond})  and
\begin{eqnarray*}
c_{n}\lim_{s \to 0}
\frac{|K|-|K_{f_{PQ},s}|}
{s^\frac{2}{n-1}}&= &D_{KL}(P_K\| Q_K)  + 2 \log \left(\frac{R}{r} \right)  \frac{as_{\pm \infty} (K) }{|K^\circ|}.
\end{eqnarray*}
If $K$ is in $C^2_+$, then this equals  $D_{KL}\big(N_{K}N_{K^{\circ}}^{-1}cm_{ K^{\circ}}\| cm_{ K}) + 2n  \log \left(\frac{R}{r} \right)$. 

\begin{eqnarray*}
c_{n}\lim_{s \to 0}
\frac{|K|-|K_{f_{QP},s}|}
{s^\frac{2}{n-1}} &=& D_{KL}(Q_K\| P_K) + 2n \log \left(\frac{R}{r} \right)  
\end{eqnarray*}
If $K$ is in $C^2_+$, then this is equal to $D_{KL}\big(N_{K}N_{K^{\circ}}^{-1}cm_{ K^{\circ}}\| cm_{ K}) + 2n  \log \left(\frac{R}{r} \right)$. 
\end{corollary}
\vskip 4mm
\noindent
{\bf Proof.} 
Note that $r=R$ iff $K$ is a Euclidean ball with radius $r$. Then the right hand sides of the identities in the corollary are equal to $0$ and $f_{PQ}$ and $f_{QP}$  are identically equal to $\infty$. Therefore, for all $s\geq 0$, $K_{f_{PQ},s} = K$ and 
$K_{f_{QP},s} = K$ and hence for all $s\geq 0$, $|K|-|K_{f_{PQ},s}| =0$ and $|K|-|K_{f_{QP},s}| =0$. Therefore, the corollary holds trivially in this case.
\par 
Assume now that $r <R$. Then
$$
1 \leq \frac{ R^{2n}|K| \ \kappa_K(x)} {r^{2n} |K^\circ| \  \langle x, N_{K}(x) \rangle^{n+1}}   \leq \left(\frac{R}{r}\right)^{4n},
$$
and we get for all $x \in \partial K$ that
$$
f_{PQ}(x) \geq 
\left(\frac { |K^\circ|  r^{n-1}}{2  \log\left(\frac{R}{r}\right)}\right)^\frac{n-1}{2} >0.
$$
Also, for all $x \in \partial K$, $ \left(\frac { |K^\circ|  r^{n-1}}{2  \log\left(\frac{R}{r}\right)}\right)^\frac{n-1}{2}  \leq M_{f_{PQ}}(x) \leq \infty$ and therefore $f_{PQ}$ satisfies the integrability condition (\ref{intcond}). 
The proof of the corollary then follows immediately from Theorem \ref{thm:surface}.
If $K$ is in $C^2_+$,  
\par
Similarly for $f_{QP}$.
\par
If $K$ is in $C^2_+$,  condition (\ref{condition}), holds. We can take 
\begin{equation} \label{rR}
r= \text{inf}_{x \in \partial K} \min_{1 \leq i \leq n-1}  r_i(x) 
\ \  \text{and}  \ \  R= \sup_{x \in \partial K} \max_{1 \leq i \leq n-1}  r_i(x), 
\end{equation}    
where for $x \in \partial K$, $r_i(x)$, $1 \leq i \leq n-1$ are the principal radii of curvature.
\vskip 5mm

For convex bodies $K$ and $ K_i$, $i=1, \cdots, n$, define
$$
\tilde{f}(N^{-1}_K(u))=f_K(u)^{\frac{n-2}{2}}[f_p(K_1,u)\cdots
f_p(K_n,u)]^{\frac{1-n}{2(n+p)}},
$$
where $ f_p(K,u)= h_K(u)^{1-p}f_K(u)$.
\vskip 3mm
\noindent

\begin{corollary}\label{int:mixed}
Let
$K$ and $ K_i$, $i=1, \cdots, n$,  be convex bodies in $C^2_+$.
Then
 \begin{eqnarray*}
&& \hskip -10mm \frac{c_n}{n \left( \prod_{i=1}^n |K_i| |K_i^\circ|^\frac{P}{n} \right) ^\frac{1}{n+p}}   \lim
_{s\rightarrow 0}  \frac{|K|-|K_{\tilde{f},s}|}{
s^{\frac{2}{n-1}}} = \\
 && \hskip 5mm  \text{Exp} \left( - \frac {n}{n+p}  D_{\frac{p}{n+p}}(P_{K_1} \times \dots \times P_{K_n} \|Q_{K_1} \times \dots \times Q_{K_n})\right),
\end{eqnarray*}
and 
 \begin{eqnarray*}
&& \hskip -10mm \frac{c_n}{n \left( \prod_{i=1}^n |K_i| |K_i^\circ|^\frac{P}{n} \right) ^\frac{1}{n+p}}   \lim
_{s\rightarrow 0}  \frac{|K|-|K_{\tilde{f},s}|}{
s^{\frac{2}{n-1}}} = \\
 && \hskip 5mm  \text{Exp} \left( - \frac {p}{n+p}  D_{\frac{n}{n+p}}(Q_{K_1} \times \dots \times Q_{K_n} \|P_{K_1} \times \dots \times P_{K_n})\right).
\end{eqnarray*}

\end{corollary}

\vskip 4mm
\noindent
{\bf Proof.} Again, the proof  follows immediately from Theorems \ref{thm:surface} and \ref{mixed-aff=renyi}.

\vskip 5mm

\noindent
{\bf Remark.}
\par
\noindent
It was shown in \cite{WY2010} that 
for  a  convex body $K$ in $\mathbb{R}^n$ with 
$C^2_+$-boundary 
 \begin{equation}\label{Form:0:0:0}
 \lim
_{s\rightarrow 0} c_n \frac{|K^{f,s}|-|K|}{
s^{\frac{2}{n-1}}}=\int_{\partial K} \frac{\kappa
_K(x)^{\frac{1}{n-1}}}{f(x)^{\frac{2}{n-1}}}\,d\mu _K(x),
\end{equation} 
where $c_n=2 |B_2^{n-1}|^{\frac{2}{n-1}}$ and  $f:
\partial K \rightarrow \mathbb{R}$ is an integrable function such  that  $f \geq c$ $\mu _K $-almost everywhere.  $c>0$ is  a constant.
Using  (\ref{Form:0:0:0}),  similar geometric interpretations of R\'enyi divergence  can be obtained via the illumination surface body instead of the surface body. We can use the same functions 
as in  Corollary \ref{cor:asp}, Corollary \ref{cor:dkl} and Corollary \ref{int:mixed}. We will also have to assume  that $K$ is in $C^2_+$.

\vskip 5mm
In \cite{PW1}, the following   new affine invariants $\Omega_K$  were introduced and its relation to the relative entropies were established:
\par
\noindent
Let $K, K_1, \dots, K_n$ be  convex bodies  in $\mathbb R^n$, all  with centroid at the origin.  Then
$$ 
\Omega_{K} = \lim_{p\rightarrow \infty} \left(\frac{as_{p}(K)}{n |K^{\circ}|}\right)^{n+p}.
$$
and   
$$ 
\Omega_{K_1, \dots K_n} = \lim_{p\rightarrow \infty} \left(\frac{as_{p}(K_1, \dots, K_n)}{as_{\infty}(K_1, \dots, K_n)}\right)^{n+p}.
$$

It was proved in \cite{PW1} that 
for  a  convex body $K$ in $\mathbb R^n$ that is $C^{2}_{+}$
\begin{equation}\label{prop3:eq1} 
D_{KL}(P_K\|Q_K) = \log{\left( \frac{|K|}{|K^\circ|} \Omega_{K}^{-\frac{1}{n}} \right) }
\end{equation}
and
\vskip 2mm
\noindent
\begin{equation}\label{prop3:eq2} 
D_{KL}(Q_K\|P_K) = \log{ \left(  \frac{|K^{\circ}|}{|K|} \Omega_{K^{\circ}}^{-\frac{1}{n}} \right) }.
\end{equation}
Note that equation (\ref{prop3:eq1}) also followed from (\ref{omega-lim}).
Similar results hold for $\Omega_{K_1, \dots K_n} $.  We now concentrate on $\Omega_{K}$.
As shown in \cite{PW1}, these invariants can also be obtained as 
 $$ 
 \Omega_{K}^\frac{1}{n}  =  \lim_{p\rightarrow 0}  \left(\frac{as_{p}(K^{\circ})}{n|K^{\circ}|}\right)^{\frac{n+p}{p}}
 $$
and thus, 
denoting by 
$
\mathcal{A}_{K} = \lim_{p\rightarrow 0} \left(\frac{as_{p}(K)}{n |K|}\right)^\frac{n+p}{p}$, 
$\Omega_K^\frac{1}{n}=\mathcal{A}_{K^\circ}$. 
This implies   e.g. that
$$
\lim_{p \rightarrow 0} \left(\frac{as_{p}(K)^n}{as_{\frac{1}{p}}(K)^\frac{1}{n}} \frac{n^\frac{1}{n} |K^\circ|^\frac{1}{n} } {n^n |K|^n}\right)^\frac{1}{p} = 1.
$$
\vskip 5mm
Geometric interpretations in terms of $L_p$-centroid bodies were given in \cite{PW1} for the new 
affine invariants $\Omega_K$. These interpretations are in the spirit of Corollaries \ref{cor:asp},  \ref{cor:dkl} and \ref{int:mixed}: As $p \rightarrow  \infty$,  appropriately chosen volume differences of  $K$ and its $L_p$-centroid bodies make  the quantity  $\Omega_K$ appear. 
\par
Again, however, with the  $L_p$-centroid bodies,  only symmetric convex bodies in $C^2_+$ could be handled and it was needed to go to a second order expansion for the volume differences. 
\par
Now,  it follows from Corollary \ref{cor:dkl}  that  there exist such interpretations for $\Omega_K$ also for non-symmetric convex bodies and under weaker smoothness assumptions than $C^2_+$.
\par
Moreover, again already a first order expansion gives such  geometric interpretations if one uses  the surface bodies or the illumination surface bodies instead of the $L_p$-centroid bodies.

\vskip 5mm 
\noindent 
\begin{corollary} \label{cor:3} 
Let $K$ be a convex body in $\mathbb{R}^n$ such that
$0$ is the center of gravity of $K$ and such  that $K$ satisfies
(\ref{condition}) of Corollary \ref{cor:dkl}.
Let $f_{PQ}: \partial K \rightarrow \mathbb{R}$ and $f_{QP}: \partial K \rightarrow \mathbb{R}$ be as in Corollary \ref{cor:dkl}. 
Then
\begin{eqnarray*}
c_{n}\lim_{s \to 0}
\frac{|K|-|K_{f_{PQ},s}|}
{s^\frac{2}{n-1}} - 2 \log \left(\frac{R}{r} \right)  \frac{as_{\pm \infty} (K) }{|K^\circ|} &=& \log \left(\frac{  |K|}{|K^{\circ}|} \Omega_K^{-\frac{1}{n}}\right) \\
&=& \log \left(  \frac{|K|}{|K^\circ|}\mathcal{A}_{K^{\circ}}^{-1} \right) 
\end{eqnarray*}
and
\begin{eqnarray*}
c_{n}\lim_{s \to 0}
\frac{|K|-|K_{f_{QP},s}|}
{s^\frac{2}{n-1}} - 2 n  \log \left(\frac{R}{r} \right)  =  \log \left(\frac{  |K^\circ|}{|K|} \Omega_{K^\circ}^{-\frac{1}{n}}\right) = \log\left( \frac{|K^\circ|}{|K|} \mathcal{A}_{K}^{-1} \right) .
\end{eqnarray*}

\end{corollary}
\vskip 4mm
\noindent
{\bf Proof.} The proof of the corollary follows immediately from Corollary  \ref{cor:dkl}, (\ref{prop3:eq1}), (\ref{prop3:eq2} ) and the definition of $\mathcal{A}_{K}$.

 \vskip 2mm \noindent Elisabeth Werner\\
{\small Department of Mathematics \ \ \ \ \ \ \ \ \ \ \ \ \ \ \ \ \ \ \ Universit\'{e} de Lille 1}\\
{\small Case Western Reserve University \ \ \ \ \ \ \ \ \ \ \ \ \ UFR de Math\'{e}matique }\\
{\small Cleveland, Ohio 44106, U. S. A. \ \ \ \ \ \ \ \ \ \ \ \ \ \ \ 59655 Villeneuve d'Ascq, France}\\
{\small \tt elisabeth.werner@case.edu}\\ \\

\end{document}